\documentclass{amsart}
\usepackage{amssymb}
\usepackage{amscd}

\usepackage[
hyperref,
nonotes
]{degt}

\def\gcd{\QOPNAME{gcd}}
\def\lcm{\QOPNAME{lcm}}

\let\CF\cf
\def\CK{\Cal K}

\let\Fermat=\Phi
\def\FC{\Fermat^\circ}
\let\rl=R

\let\ln=L

\let\dv=V
\def\two{_{(2)}}
\def\two{_{=}}
\def\tFermat{\Fermat'}
\def\tdv{\dv'}
\def\bm{\bold m}
\def\bm{\frak m}

\let\kk=\Bbbk
\def\FF{\Bbb F}
\def\sset{\bold S}
\def\Kset{\bold K}
\def\Tset{\bold T}
\def\SS{\sset\ofD}
\def\KK{\Kset\ofD}
\def\TT{\Tset\ofD}

\def\SF#1{\SS{\dv\of{#1}}}
\def\KF#1{\KK{\dv\of{#1}}}
\def\TF#1{\TT{\dv\of{#1}}}
\def\SF#1{\sset\of{#1}}
\def\KF#1{\Kset\of{#1}}
\def\TF#1{\Tset\of{#1}}
\def\tS#1{\tilde\sset\of{#1}}
\def\ofD#1{\<#1\>}
\def\GG{\Bbb G}
\def\Pic{\QOPNAME{Pic}}

\def\of#1{(#1)}
\def\of#1{[#1]}
\def\ofm{\of{m}}

\let\Am=A
\def\Am{\mathrm{A}}

\let\Bm=B
\def\Bm{\mathrm{B}}
\def\Km{\mathrm{K}}

\def\tBm{\tilde\Bm}

\let\=\B

\let\vg=v
\let\hg=h
\let\ag=a
\let\cg=c
\let\ug=u
\let\g=g
\let\3\tilde

\let\GL\Lambda

\def\tGL{\reduced\GL}
\def\tGL{\mathrm{R}}

\def\height{\QOPNAME{ht}}

\let\inc\iota

\let\cp=\Gs      
\def\bcp{\bar\cp}
\let\cf=\psi     
\def\cpp{\tilde\varphi} 
\let\Euler\phi   
\let\co=o        
\let\cd=\Gr      
\def\Gda{\Gd\of\Ga}
\def\inf{f}      

\def\CW{{\sl CW}}
\def\diag{\QOPNAME{diag}}
\let\GM=M
\let\sm=T

\def\mat[ [#1,#2,#3], [#4,#5,#6], [#7,#8,#9] ]{\left[\begin{array}{rrr}
#1&#2&#3\\#4&#5&#6\\#7&#8&#9
\end{array}\right]}

\title
 {On the Picard group of a Delsarte surface}

\author{Alex Degtyarev}

\address{%
Department of Mathematics\\
Bilkent University\\
06800 Ankara, Turkey}

\email{degt@fen.bilkent.edu.tr}



\begin{document}

\begin{abstract}
We suggest an algorithm computing, in some cases,
an explicit generating set for the
N\'{e}ron--Severi lattice of a Delsarte surface.
\end{abstract}

\maketitle

\section{Introduction}

Throughout the paper, all algebraic varieties are over~$\C$.

\subsection{Statement of the problem}\label{s.statement}
A \emph{Delsarte surface} is a surface $\Fermat_A\subset\Cp3$ given
by a four-term equation of the form
\[\label{eq.Delsarte}
\sum_{i=0}^3\prod_{j=0}^3z_j^{a_{ij}}=0,
\]
see~\cite{Delsarte,Shioda:Delsarte}.
The restrictions to the matrix~$A:=[a_{ij}]$ are listed in
\autoref{s.Delsarte}\iref{A.integer}--\iref{A.det}.

We are interested in certain birational invariants of
Delsarte surfaces. For this reason,
we silently replace~$\Fermat_A$ with its
resolution of singularities. The particular choice of the
resolution is not important; \eg, one can take the minimal one.

For an alternative description of Delsarte surfaces, introduce the
multiplicative abelian group $\GG\cong\Z^3$ with a distinguished
generating set $t_0,t_1,t_2,t_3$ subject to the only relation
$t_0t_1t_2t_3=1$. Then, each epimorphism $\Ga\:\GG\onto G$
to a finite group~$G$
gives rise to a Delsarte surface $\Fermat\of\Ga$, see \autoref{s.Delsarte}
and \autoref{def.Delsarte}.
By an abuse of the language, an epimorphism~$\Ga$ as above
is referred to as a \emph{finite quotient} of~$\GG$.

\definition\label{def.special}
In the examples,
we will consider the following four special classes of
Delsarte surfaces, corresponding to special finite quotients
$\Ga\:\GG\onto G$\rom:
\roster
\item\label{def.Fermat}
\emph{Fermat surfaces} $\Fermat\ofm$, where an integer
$m\in\N_+$ is regarded
as the quotient projection $m\:\GG\onto\GG/m\GG$\rom;
\item\label{def.unramified}
\emph{unramified \rom(at~$\infty$\rom) Delsarte surfaces} $\Fermat\of\Ga$,
\ie, such that
$\Ga(t_0)=1$\rom;
\item\label{def.cyclic}
\emph{cyclic Delsarte surfaces} $\Fermat\of\Ga$, \ie, such that $G$
is a cyclic group\rom;
\item\label{def.diagonal}
\emph{diagonal Delsarte surfaces} $\Fermat\of{\bm}$,
where a vector $\bm:=(m_1,m_2,m_3)\in\N_+^3$ is regarded
as the quotient projection
$\bm\:\GG\onto\GG/(t_1^{m_1}=t_2^{m_2}=t_3^{m_3}=1)$.
\endroster
(To avoid the common confusion, we use $\N_+$ for the set of \emph{positive}
integers.)
Note that, in items~\iref{def.unramified} and~\iref{def.diagonal}, the definition
depends on the order of the indices, and we relate a surface $\Fermat\of\Ga$
to the corresponding class whenever it satisfies the condition after a
possible
permutation of the indices $(0,1,2,3)$.
\enddefinition

By Poincar\'{e} duality,
the N\'{e}ron--Severi lattice
$\NS(\Fermat\of\Ga)$
can be regarded
as a subgroup of the homology group $H_2(\Fermat\of\Ga)/\!\Tors$.
Our primary interest is the extent to which $\NS(\Fermat\of\Ga)$ is generated by
the components of a certain `obvious' divisor $\dv\of\Ga\subset\Fermat\of\Ga$,
see \autoref{s.dv}.
(In the case of Fermat surfaces,
this divisor $\dv$ is essentially constituted by the
lines contained in the surface.)
To this end, we consider the inclusion homomorphism
$\inc_*\:H_2(\dv\of\Ga)\to\NS(\Fermat\of\Ga)$ and introduce the
groups
\[
\SF\Ga:=\Im\inc_*,\quad
\KF\Ga:=\Ker\inc_*,\quad
\TF\Ga:=\Tors\bigl(\NS(\Fermat\of\Ga)/\SF\Ga\bigr).
\label{eq.of.alpha}
\]
We compute the two latter groups, which are birational invariants of the
surface.

The motivation for our interest is Shioda's algorithm~\cite{Shioda:Delsarte}
computing the Picard rank
$\rho(\Fermat\of\Ga)$.
In some cases (most notably, if $\ls|G|$ is prime to~$6$,
\CF. \autoref{cor.generated} below), this computation implies that
$\NS(\Fermat\of\Ga)\otimes\Q=\SF\Ga\otimes\Q$,
\ie, $\NS(\Fermat\of\Ga)$ is generated by the components of~$\dv\of\Ga$
over~$\Q$; hence, a natural question is if this generation property still
holds \emph{over the integers}, \ie, if $\TF\Ga=0$.
We answer this question in the affirmative for
a few special classes of
surfaces, while showing that, in general, the answer is in the negative.

\subsection{Principal results}\label{s.results}
Introduce the following subgroups of~$\GG$:
\roster*
\item
$\GG_{ij}$ is generated by $t_i$ and~$t_j$, $i,j=0,1,2,3$;
\item
$\GG_i$ is generated by $t_it_j$ and $t_it_k$, $i=1,2,3$ and
$\{i,j,k\}=\{1,2,3\}$;
\item
$\GG\two:=\sum_i\GG_i$ is generated by $t_1t_2$, $t_1t_3$, and $t_2t_3$.
\endroster
Given a finite quotient $\Ga\:\GG\onto G$, denote $G_*:=G/\Ga(\GG_*)$ (for a
subscript $*$ of the form $ij$, $i$, or~$=$) and let
$\Gd\of\Ga:=\ls|G\two|-1\in\{0,1\}$.
(In more symmetric terms,
$\GG_i$ depends only on the partition $\{0,i\}\cup\{j,k\}$ of the index
set, and $\GG\two$ is generated by all products $t_it_j$, $i,j=0,1,2,3$;
one has $[\GG:\GG\two]=2$.)

Recall that the \emph{length} $\ell(A)$ of a finitely generated
abelian group~$A$ is the
minimal number of generators of~$A$, and the \emph{exponent} $\exp A$ of a
finite abelian group~$A$ is the minimal positive integer~$m$ such that
$mA=0$.
For a finite quotient $\Ga\:\GG\onto G$,
the exponent $\exp G$ is the minimal positive integer~$m$ such that
$m\GG\subset\Ker\Ga$, and
we can also define the \emph{height}
$\height\Ga:=\exp G/n$, where $n$ is the maximal integer such that
$\Ker\Ga\subset n\GG$.
Note that $(\exp G)^3\!/\ls|G|$ is an integer dividing $(\height\Ga)^2$.

The principal results of the paper (combined with those of~\cite{degt:Fermat})
are stated below, with references to the proofs given in the statements.

\theorem[see \cite{degt:Fermat} and \autoref{proof.pi1}]\label{th.pi1}
For any finite quotient $\Ga\:\GG\onto G$,
one has
\[*
\pi_1(\Fermat\of\Ga)=H_1(\Fermat\of\Ga)=\Ker\Ga\big/\prod(\GG_{ij}\cap\Ker\Ga),
\]
the product running over all pairs $0\le i<j\le3$.
This group is trivial for any of the four special classes of Delsarte
surfaces introduced in \autoref{def.special}.
In general, the group $\pi_1(\Fermat\of\Ga)$ is cyclic
and its order $\ls|\pi_1(\Fermat\of\Ga)|$ divides $\height\Ga$.
\endtheorem

\theorem[see \autoref{proof.main}]\label{th.main}
For any finite quotient $\Ga\:\GG\onto G$,
one has
\[*
\rank\KF\Ga=
 \sum_{0\le i<j\le3}\ls|G_{ij}|+\sum_{1\le i\le3}\ls|G_{i}|
 -3-\Gd\of\Ga.
\]
Besides, one has
$\ell(\TF\Ga)\le6+\Gd\of\Ga$ and
$\exp\TF\Ga$ divides $(\exp G)^3\!/\ls|G|$.
\endtheorem

\addendum[see \autoref{proof.lattice}]\label{ad.lattice}
As a lattice, $\SF\Ga=H_2(\dv\of\Ga)/\ker$, where
$\ker$ is the kernel $\ker H_2(\dv\of\Ga)$ of the intersection index form.
\endaddendum

Note that \autoref{th.main} is merely an estimate on the size of
the torsion $\TF\Ga$, most
interesting
being the fact that the length
of this group is
universally bounded.
A better estimate is found in \autoref{lem.torsion},
and a precise, although not very efficient, algorithm
for computing this group is given by~\eqref{eq.Tors}.
A few examples, showing the sharpness of most estimates, are considered in
\autoref{S.examples}.
It appears that
there should be better bounds taking into account
the group $\pi_1(\Fermat\of\Ga)$, see \autoref{rem.pi1}.

Note also that the rank formula in \autoref{th.main} states, essentially,
that the rank $\rank\KF\Ga$ is `minimal possible'. More precisely, $G$ acts
on $\Fermat\of\Ga$ and $\dv\of\Ga$, see \autoref{s.Delsarte},
and the
space $H_2(\Fermat\of\Ga;\C)$ splits into
multi-eigenspaces, which are all of dimension at most~$1$,
see~\cite{Shioda:Fermat,Shioda:Delsarte}.
Comparing the dimensions (or using the explicit description of the kernel,
see \autoref{proof.main}), one can see that each eigenspace present in
$H_2(\dv\of\Ga;\C)$ is mapped \emph{epimorphically} onto the corresponding
eigenspace in $H_2(\Fermat\of\Ga;\C)$.

\theorem\label{th.free}
One has $\TF\Ga=0$ in each of the following three cases\rom:
\roster
\item\label{free.Fermat}
Fermat surfaces~$\Fermat\of\Ga$, $\Ga=m\in\N_+$,
see~\cite{degt:Fermat} or
\rom{\autoref{proof.Fermat}}\rom;
\item\label{free.unramified}
Delsarte surfaces unramified at~$\infty$,
see~\cite{degt:Fermat}\rom;
\item\label{free.cyclic}
cyclic Delsarte surfaces, see \rom{\autoref{S.cyclic}}.
\endroster
Besides, one has the following stronger bound\rom:
\roster[4]
\item\label{free.diagonal}
if
$\Fermat\of\Ga$ is
diagonal,
$\Ga=(m_1,m_2,m_3)\in\N_+^3$,
then $\ell(\TF\Ga)\le\Gd\of\Ga$
and the order $\ls|\TF\Ga|$ divides
$\lcm_{1\le i<j\le3}(\gcd(m_i,m_j))/\!\gcd(m_1,m_2,m_3)$,
see \rom{\autoref{proof.diagonal}}.
\endroster
\endtheorem

For Fermat surfaces, the primitivity statement was suggested
in~\cite{Shioda:Fermat,Shioda:Fermat1}, and it was verified numerically
in~\cite{Shioda:Fermat2} for all values of~$m$ prime to~$6$ in the range
$5\le m\le100$.
For cyclic Delsarte surfaces, \autoref{th.free}\iref{free.cyclic}, the
statement was conjectured in~\cite{Shimada.Takahashi:primitivity}, where it
was verified for all cyclic quotients $\Ga\:\GG\onto G$ with
$\ls|G|\le50$.

\corollary[see \autoref{proof.generated}]\label{cor.generated}
Let $\Ga\:\GG\onto G$ be a finite quotient
with $\Fermat\of\Ga$
in one of the four special classes introduced in \autoref{def.special},
and
assume that
$\ls|G|$ is prime to~$6$.
Then the group
$\Pic\Fermat\of\Ga=\NS(\Fermat\of\Ga)$ is generated by the components of
$\dv\of\Ga$.
In other words, $\NS(\Fermat\of\Ga)=H_2(\dv\of\Ga)/\ker$ as a lattice,
see \autoref{ad.lattice}.
\endcorollary

It is worth emphasizing that, since both the action of~$G$ (obvious) and the
intersection matrix of $\dv\of\Ga$ (see, \eg,
\cite{Shioda:Fermat2}) are known,
\autoref{cor.generated} gives us a complete description of the
N\'{e}ron--Severi group $\NS(\Fermat\of\Ga)$, including the lattice structure
and the action of $G\subset\Aut\Fermat\of\Ga$.
In general, if $\NS(\Fermat\of\Ga)\otimes\Q=\SF\Ga\otimes\Q$ but
$\TF\Ga\ne0$, the lattice structure can be recovered using the algorithm
outlined in \autoref{s.torsion}.

\subsection{Contents of the paper}
In \autoref{S.prelim},
we introduce Delsarte surfaces~$\Fermat$ and the `obvious' divisors
$\dv\subset\Fermat$ and discuss their description in terms of ramified
coverings of the plane.
In \autoref{S.topology}, most principal results of the paper are reduced to
the problem of analyzing the integral torsion of a certain
Alexander module, see~\eqref{eq.Tors} and~\eqref{eq.rank}.
Most
result are proved in \autoref{S.main}; an exception is the
case of cyclic Delsarte surfaces, which is treated separately (and slightly
differently) in \autoref{S.cyclic}.
Finally, in \autoref{S.examples}, we discuss a few numeric examples (obtained
from experiments with small random matrices), illustrating the sharpness of
most bounds on the one hand and the complexity of the general problem on
the other.

\section{Preliminaries}\label{S.prelim}

\subsection{Conventions}\label{s.notation}
The notation
$\Tors A$ stands for the $\Z$-torsion of an abelian group~$A$.
We emphasize that $\Tors$ always refers to the \emph{integral} torsion,
even if $A$ is a module over a larger ring.
This convention applies also to the rank $\rank A$ and length $\ell(A)$:
we regard~$A$ as an abelian group.
We abbreviate $A/\!\Tors:=A/\!\Tors A$.

We denote by~$\cf_m(t)$ the cyclotomic polynomial of order~$m$, \ie,
the irreducible (over~$\Q$) factor of $t^m-1$ that does not divide $t^n-1$
for $1\le n<m$. We also make use of the polynomials
$\cpp_m(t):=(t^m-1)/(t-1)$, $m\in\N_+$.

Unless stated otherwise, all homology and cohomology groups have coefficients
in~$\Z$. Since all spaces involved have homotopy type of \CW-complexes,
the
choice of a theory is not important; for example, one can use
singular (co-)homology.

Given a closed oriented $4$-manifold~$X$, we
identify $H^2(X)=H_2(X)$
by means of Poincar\'{e} duality.
In particular, if $X$ is a smooth compact complex analytic surface, we regard
the N\'{e}ron--Severy lattice $\NS(X)$ as a sublattice of $H_2(X)/\!\Tors$
(with the usual intersection index pairing), so that
a divisor $D\subset X$ is represented by its (topological) fundamental class
$[D]\in H_2(X)/\!\Tors$.

Given a smooth compact analytic surface~$X$ and a divisor $D\subset X$,
we denote by $\SS{D}\subset\NS(X)$
the subgroup generated by the
irreducible components of~$D$. In other words,
\[*
\SS{D}=\Im[\inc_*\:H_2(D)\to H_2(X)/\!\Tors],
\]
where $\inc\:D\into X$ is the inclusion.
We will also consider the groups
\[*
\TT{D}:=\Tors(\NS(X)/\SS{D}),\quad
\KK{D}:=\Ker[\inc_*\:H_2(D)\to H_2(X)/\!\Tors],
\]
which are birational invariants of the pair $(X,D)$.
More precisely, if $\Gs\:X'\to X$ is a blow-down map and $D':=\Gs^*D$,
then $\Gs^*$ and $\Gs_*$ induce isomorphisms
\[
\NS(X')/\SS{D'}=\NS(X)/\SS{D},\quad
\TT{D'}=\TT{D},\quad
\KK{D'}=\KK{D}.
\label{eq.sigma}
\]

\subsection{Delsarte surfaces\noaux{ (see~\cite{Shioda:Delsarte})}}\label{s.Delsarte}
Consider the surface~$\Fermat_A$ given by~\eqref{eq.Delsarte}, where
the exponent matrix $A:=[a_{ij}]$ is assumed to
satisfy the following conditions:
\roster
\item\label{A.integer}
each entry~$a_{ij}$, $0\le i,j\le3$, is a non-negative integer;
\item\label{A.0}
each column of~$A$ has at least one zero;
\item\label{A.const}
$(1,1,1,1)^t$ is an eigenvector of~$A$,
\ie, $\sum_{j=0}^3a_{ij}=\Gl=\const(i)$;
\item\label{A.det}
$A$ is non-degenerate, \ie, $\det A\ne0$.
\endroster
Condition~\iref{A.0} asserts that the surface does not contain a coordinate
plane, and \iref{A.const} makes~\eqref{eq.Delsarte} homogeneous, the degree
being the eigenvalue~$\Gl$.

Following~\cite{Shioda:Delsarte}, introduce the cofactor matrix
$A^*:=(\det A)A\1$ and let
\[*
d:=\gcd(a_{ij}^*),\quad
m:=\ls|\det A|/d,\quad
B:=mA\1=\pm d\1A^*.
\]
Denoting by~$\Fermat\ofm$ the Fermat surface $\{z_0^m+\ldots+z_3^m=0\}$, we
have maps
\[*
\Fermat\ofm\overset{\pi_B}\longto\Fermat_A
\overset{\pi_A}\longto\Fermat:=\Fermat\of1
\]
given by
\[*
\pi_B\:(z_i)\mapsto\Biggl(\prod_{j=0}^3z_j^{b_{ij}}\Biggr),\quad
\pi_A\:(z_i)\mapsto\Biggl(\prod_{j=0}^3z_j^{a_{ij}}\Biggr).
\]
Both maps are ramified coverings;
$\pi_A$ and $\pi_B\circ\pi_A\:(z_i)\mapsto(z_i^m)$
are ramified over the union
$\rl:=\rl_0+\rl_1+\rl_2+\rl_3\subset\Fermat$ of the
traces of the coordinate
planes, $\rl_i:=\Fermat\cap\{z_i=0\}$.

The fundamental group $\pi_1(\Fermat\sminus\rl)$ is abelian and,
by Poincar\'{e}--Lefschetz duality,
there are canonical isomorphisms
\[*
\pi_1(\Fermat\sminus\rl)=H^2(\rl)/H^2(\Fermat)=\GG,
\]
where
$\GG$ is the abelian group introduced in \autoref{s.statement} and
a generator $t_i\in\GG$ evaluates to the
Kronecker symbol~$\Gd_{ij}$ on the fundamental class $[R_j]$
(with its canonical complex orientation).
Thus, away from the ramification locus~$\rl$, the
unramified topological covering $\pi_A$ is uniquely
determined by a finite index subgroup of~$\GG$,
\viz. the image of $\pi_1(\Fermat_A\sminus\pi_A\1(\rl))$,
or, equivalently, by a finite
quotient $\Ga\:\GG\onto G$.

Due to condition~\iref{A.const} above,
$A$ and~$B$
can be regarded as endomorphisms of~$\GG$, inducing endomorphisms
$A_m,B_m\:\GG/m\GG\to\GG/m\GG$. Obviously, one has
\[*
m\GG\subset\Ker\Ga,\quad
\MG:=\Ker\Ga/m\GG=\Ker B_m=\Im A_m,\quad
\Im B_m=\Ker A_m,
\]
and $\Fermat_A$ is birationally isomorphic to $\Fermat\ofm/\MG$, where
a generator $t_i\in\GG/m\GG$ acts
on $\Fermat\ofm$ by multiplying the $i$-th coordinate by a fixed primitive
$m$-th root of unity.

Summarizing, we can disregard the original exponent matrix~$A$ and
equation~\eqref{eq.Delsarte} and adopt the following definition,
\CF.~\cite{Shimada.Takahashi:primitivity}.

\definition\label{def.Delsarte}
Given a finite
quotient $\Ga\:\GG\onto G$, the \emph{Delsarte surface} $\Fermat\of\Ga$ is
defined as
(any) smooth analytic compactification of the
(unramified) covering of
the
complement
$\Fermat\sminus\rl$ corresponding to~$\Ga$.
\enddefinition

Since
the invariants that we are interested in are of a birational nature,
\CF. \eqref{eq.sigma}, the particular choice of the compactification
in \autoref{def.Delsarte} is not important.
It is fairly obvious that any surface $\Fermat\of\Ga$ is a resolution of
singularities of the projective surface given by an appropriate
equation~\eqref{eq.Delsarte}; however, we do not use this fact.
For the covering Fermat surface~$\Fermat\ofm$,
we can merely take $m=\exp G$ or any
multiple thereof, so that $m\GG\subset\Ker\Ga$.

\subsection{The divisor \pdfstr{V[a]}{$\dv\of\Ga$}}\label{s.dv}
Fix a finite quotient $\Ga\:\GG\onto G$ and let
$\pi\:\Fermat\of\Ga\to\Fermat$ be the covering projection.
Consider the lines $\ln_i:=\Fermat\cap\{z_0+z_i=0\}$, $i=1,2,3$,
let $\ln:=\ln_1+\ln_2+\ln_3$,
and define the divisors
\[*
\rl_*\of\Ga:=\pi^*\rl_*,\quad
\ln_*\of\Ga:=\pi^*\ln_*,\quad
\dv\of\Ga:=\rl\of\Ga+\ln\of\Ga
\]
on~$\Fermat\of\Ga$. (Here, the subscript $_*$ is either empty or an
appropriate index in the range $0,\ldots,3$.)
To avoid excessive nested parentheses, introduce the shortcuts
\[*
\SF\Ga:=\SS{\dv\of\Ga},\quad
\KF\Ga:=\KK{\dv\of\Ga},\quad
\TF\Ga:=\TT{\dv\of\Ga},
\]
\CF.~\eqref{eq.of.alpha}, and let
$\FC\of\Ga:=\Fermat\of\Ga\sminus\dv\of\Ga$.
We recall that the pull-back of each~$\ln_i$, $i=1,2,3$, in the covering Fermat
surface~$\Fermat\ofm$ splits into $m^2$ `obvious' straight lines,
\viz.
\[
\aligned
\ln_1(\zeta,\eta):&\quad(r:\Go\zeta r:s:\Go\eta s),\\
\ln_2(\zeta,\eta):&\quad(r:s:\Go\zeta r:\Go\eta s),\\
\ln_3(\zeta,\eta):&\quad(r:s:\Go\eta s:\Go\zeta r),
\endaligned
\label{eq.lines}
\]
where $(\zeta,\eta)$ is a pair of $m$-th roots of unity (parametrizing the
$m^2$ lines within each of the three families), $\Go:=\exp(\pi i/m)$ is an
$m$-th root of~$-1$, and $(r:s)$ is a point in~$\Cp1$,
\CF.~\cite{Shioda:Fermat2}.
Thus, the components of $\dv\of\Ga$ are the images of the $3m^2$ straight
lines contained in the covering Fermat surface~$\Fermat\ofm$,
the components of the
ramification locus of the covering $\Fermat\of\Ga\to\Fermat$, and the
exceptional divisors arising from the resolution of singularities.


\section{The topology of a Delsarte surface}\label{S.topology}

In this section, we discuss a few simplest topological properties of the
Delsarte surface $\Fermat\of\Ga$ and divisor $\dv\of\Ga$ defined by a finite
quotient $\Ga\:\GG\onto G$.
In particular, we reduce most statements to the study of certain modules
$\Am\of\Ga$ or $\Bm\of\Ga$.

\subsection{The fundamental group: proof of \autoref{th.pi1}}\label{proof.pi1}
The expression for the group $\pi_1(\Fermat\of\Ga)$ in terms of~$\Ga$ is found
in~\cite{degt:Fermat}, and the statement that
$\pi_1(\Fermat\of\Ga)=0$
for Fermat surfaces and unramified or diagonal Delsarte surfaces is
immediate. We postpone the case of cyclic Delsarte surfaces till
\autoref{proof.pi1.cyclic},
where the necessary framework is introduced.

In general, we can assume that
the kernel
$\Ker\Ga$ is generated by three vectors
$v_i:=t_1^{m_{i1}}t_2^{m_{i2}}t_3^{m_{i3}}$, $i=1,2,3$, so that the matrix
$[m_{ij}]$ is upper triangular,
\[*
[m_{ij}]=\bmatrix
m_{11}&m_{12}&m_{13}\\
0&m_{22}&m_{23}\\
0&0&m_{33}
\endbmatrix.
\]
Then $\GG_{23}\cap\Ker\Ga$ contains $v_3$ and~$v_2$, and
$\GG_{13}\cap\Ker\Ga$ contains
$v_3$ and a product of the form $v_1^rv_2^s$, $r\ne0$.
Hence, $\pi_1(\Fermat\of\Ga)$ is a cyclic group (generated by~$t_1$) of order
at most~$r$.
On the other hand,
from the expression in the statement,
it is clear that $\pi_1(\Fermat\of\Ga)$
is a subquotient of the group $n\GG/m\GG$
of exponent $\height\Ga$,
where $m:=\exp G$
and $n$ is as in the
definition of $\height\Ga$, see \autoref{s.results}.
\qed

\subsection{The reduction}\label{s.reduction}
Our proof of Theorems~\ref{th.main} and~\ref{th.free}
is based on the following homological
reduction of the problem.

\theorem\label{th.reduction}
Let $D$ be a divisor in a smooth compact analytic surface~$X$, and let
$K(X,D):=\Ker[\kappa_*\:H_1(X\sminus D)\to H_1(X)]$ be the kernel of the
homomorphism~$\kappa_*$
induced by the inclusion.
Then there are
canonical
isomorphisms
\[*
\Tors K(X,D)=\Hom(\TT{D},\Q/\Z),\qquad
K(X,D)/\!\Tors=\Hom(\KK{D},\Z).
\]
\endtheorem

\proof
The inclusion homomorphism $\kappa_*\:H_1(X\sminus D)\to H_1(X)$
is Poincar\'{e} dual to
the homomorphism~$\Gb$ in the following exact sequence of pair $(X,D)$:
\[*
\longto H^2(X)\overset{\inc^*}\longto H^2(D)\longto
H^3(X,D)\overset\Gb\longto H^3(X)\longto.
\]
Hence, $K(X,D)=\Coker\inc^*$, and both statements are immediate,
\CF.~\cite{degt:Fermat}, using the definition of the $\Ext$ groups
in terms of projective resolutions and the
canonical isomorphism $\Ext(A,\Z)=\Hom(A,\Q/\Z)$ for any finite abelian
group~$A$.
\endproof

\subsection{The modules \pdfstr{A[a]}{$\Am\of\Ga$} and \pdfstr{B[a]}{$\Bm\of\Ga$}}\label{s.modules}
The groups $H_1(\FC\of\Ga)=H_1(\Fermat\of\Ga\sminus\dv\of\Ga)$
for Delsarte surfaces were
computed in~\cite{degt:Fermat}, using the covering $\FC\of\Ga\to\FC$.
Let
\[*
\GL:=\Z[\GG]=\Z[t_1^{\pm1},t_2^{\pm1},t_3^{\pm1}]=
 \Z[t_0,t_1,t_2,t_3]/(t_0t_1t_2t_3-1).
\]
be the ring of Laurent polynomials, and consider the homomorphism
$\partial\:\Am\of0\to\GL$ of $\GL$-modules defined as follows:
$\Am\of0$ is the $\GL$-module generated by six elements $\ag_i$, $\cg_j$,
$i,j=1,2,3$, subject to the relations
\begin{gather}
(t_2t_3-1)\cg_1=(t_1t_3-1)\cg_2=(t_1t_2-1)\cg_3=0,\label{rel.bb}\\
 (t_3-1)\cg_1+(t_3-1)\ag_2-(t_2-1)\ag_3=0,\label{rel.a1}\\
 (t_3-1)\cg_2+(t_3-1)\ag_1-(t_1-1)\ag_3=0,\label{rel.a2}\\
 (t_1-1)\cg_3+(t_1-1)\ag_2-(t_2-1)\ag_1=0,\label{rel.a3}
\end{gather}
and $\partial$ is
\[
\partial\ag_i=(t_i-1),\quad
\partial\cg_j=0,\quad i,j=1,2,3.
\label{eq.d1}
\]
For an epimorphism $\Ga\:\GG\onto G$, let $\GL\of\Ga:=\Z[G]$.
The induced ring homomorphism
$\GL\onto\GL\of\Ga$ makes $\GL\of\Ga$ a $\GL$-module, and we define
$\Am\of\Ga:=\Am\of0\otimes_\GL\GL\of\Ga$.
In other words, $\Am\of\Ga$ is obtained from $\Am\of0$ by adding
to~\eqref{rel.bb}--\eqref{rel.a3} the defining relations of~$G$ in the basis
$\{t_1,t_2,t_3\}$.
Then, the computation in~\cite{degt:Fermat} can be summarized in the form of
an exact sequence
\[*
0\longto
H_1(\FC\of\Ga)\longto
\Am\of\Ga\overset{\partial}\longto
\GL\of\Ga\longto\Z\longto0.
\]
The
homomorphism~$\kappa_*$ in \autoref{th.reduction} factors through the free
abelian group
\[*
H_1(\Fermat\of\Ga\sminus\rl\of\Ga)=\pi_1(\Fermat\of\Ga\sminus\rl\of\Ga)
 =\Ker\Ga\cong\Z^3.
\]
The homology~$H_0$ and~$H_1$ of
the space $\Fermat\of\Ga\sminus\rl\of\Ga$ are computed
by the complex
$0\to\Am\of\Ga/\Bm\of\Ga\to\GL\of\Ga\to0$,
where $\Bm\of\Ga\subset\Am\of\Ga$ is the $\GL\of\Ga$-submodule generated
by $\cg_1,\cg_2,\cg_3$.
Summarizing, we can restate \autoref{th.reduction} as follows:
\begin{gather}
\Hom(\TF\Ga,\Q/\Z)=\Tors H_1(\FC\of\Ga)=\Tors\Am\of\Ga=\Tors\Bm\of\Ga,
 \label{eq.Tors}\\
\rank\KF{\Ga}=\rank\Am\of\Ga-\ls|G|+1=\rank\Bm\of\Ga+3.
 \label{eq.rank}
\end{gather}

\subsection{Generators of the torsion}\label{s.torsion}
An explicit generating set for the primitive hull
$\tS\Ga:=(\SF\Ga\otimes\Q)\cap\NS(\Fermat\of\Ga)$ can be described in terms
of the discriminant form.
We outline this description, in the hope that it may be useful in the
future.

The lattice $\SF\Ga$ has a vector of positive square (\eg, the hyperplane
section class); hence, the Hodge index theorem implies
that $\SF\Ga$ is non-degenerate and its dual group
$\sset^*$ can be identified with a subgroup of $\SF\Ga\otimes\Q$:
\[*
\sset^*:=\Hom(\SF\Ga,\Z)=\{x\in\SF\Ga\otimes\Q\,|\,
 \text{$x\cdot y\in\Z$ for all $y\in\SF\Ga$}\}.
\]
This identification gives rise to an inclusion $\SF\Ga\subset\sset^*$ and
to the \emph{discriminant group} $\discr\SF\Ga:=\sset^*\!/\SF\Ga$,
see~\cite{Nikulin:forms}.
The latter is a finite abelian group equipped with a
non-degenerate symmetric $\Q/\Z$-valued
bilinear form, \viz. the descent of the $\Q$-valued extension of the
intersection index form from $\SF\Ga$ to~$\sset^*$.
Since $\tS\Ga$ is also an integral lattice,
there are natural inclusions
\[*
\SF\Ga\subset\tS\Ga\subset\tilde\sset^*:=\Hom(\tS\Ga,\Z)\subset\sset^*;
\]
hence,
the extension $\tS\Ga\supset\SF\Ga$ is
uniquely determined by
either of the subgroups
\[*
\CK:=\tS\Ga/\SF\Ga\subset
 \CK^\perp:=\tilde\sset^*\!/\SF\Ga\subset
 \discr\SF\Ga.
\]
Indeed,
the subgroups $\CK\subset\CK^\perp$ are the orthogonal complements of each
other (in particular, $\CK$ is isotropic), and
\[*
\tS\Ga=\{x\in\SF\Ga\otimes\Q\,|\,x\bmod\SF\Ga\in\CK\}.
\]
For further details concerning discriminant forms and lattice extensions,
see~\cite{Nikulin:forms}.

Consider the $\GL\of\Ga$-module $\tBm\of\Ga$ generated by $\cg_1,\cg_2,\cg_3$
subject to relations~\eqref{rel.bb}.
The geometric description found in~\cite{degt:Fermat} establishes
a canonical, up to the coordinate action of~$\GG$, homomorphism
$\tBm\of\Ga\to H^2(\dv\of\Ga)$ of $\GL\of\Ga$-modules, which restricts to an
isomorphism $\tBm\of\Ga=H^2(\ln'\of\Ga)$, where $\ln'\of\Ga$ is the
\emph{proper} transform of~$\ln$ in~$\Fermat\of\Ga$.
If $\Ga=m\in\N_+$, the reference point in~$\Fermat\ofm$ can be chosen so that
\[
\cg_1\mapsto[\ln_1(1,\Go^{-2})]^*,\quad
\cg_2\mapsto[\ln_2(1,\Go^{-2})]^*,\quad
\cg_3\mapsto[\ln_3(1,1)]^*,
\label{eq.hom}
\]
see~\eqref{eq.lines} for the notation; in general, we use, in addition,
the natural identifications
$\tBm\of\Ga=\tBm\ofm\otimes_\GL\GL\of\Ga$ and
$H^2(\ln'\of\Ga)=H^2(\ln\ofm)\otimes_\GL\GL\of\Ga$.

Consider the modules
\[*
\Km':=\Ker[\tBm\of\Ga\to\Bm\of\Ga]\subset
\Km:=\Ker[\tBm\of\Ga\to\Bm\of\Ga/\!\Tors].
\]
It is immediate from the construction (with~\eqref{eq.Tors} taken into
account) that the group $\Km/\Km'$ is canonically isomorphic to
$\sset^*\!/\tilde\sset^*$. The homomorphism $\Km\to\discr\SF\Ga$ is
easily computed using~\eqref{eq.hom} and the intersection matrix of the
components of~$\dv\of\Ga$, see, \eg,~\cite{Shioda:Fermat2}, and the subgroup
$\CK^\perp\subset\discr\SF\Ga$ defining the extension
$\tS\Ga\supset\SF\Ga$ as described above is found as the image of~$\Km'$.

\section{Proof of \autoref{th.main}}\label{S.main}

Throughout this section, we consider a finite quotient
$\Ga\:\GG\onto G$ and fix the notation $m:=\exp G$.

\subsection{Alternative proof of \autoref{th.free}\iref{free.Fermat}}\label{proof.Fermat}
This proof repeats almost literally the one found in~\cite{degt:Fermat},
except that we analyze the module~$\Bm\of\Ga$ instead of~$\Am\of\Ga$. This
analysis (slightly more thorough than in~\cite{degt:Fermat})
is used in the sequel.

Assume that $\Ga=m\:\GG\onto G=\GG/m\GG$
and consider the filtration
\[*
0=\Bm_0\subset\Bm_1\subset\Bm_2\subset\Bm_3\subset\Bm_4:=\Bm\of\Ga,
\]
where
\roster*
\item
$\Bm_3$ is generated by
$\cg_1':=(t_3-1)\cg_1$, $\cg_2':=(t_3-1)\cg_2$, $\cg_3':=(t_1-1)\cg_3$,
\item
$\Bm_2$ is generated by
$\cg_1'':=(t_1-1)\cg_1'$, $\cg_2'':=(t_2-1)\cg_2'$, $\cg_3'':=(t_3-1)\cg_3'$,
and
\item
$\Bm_1$ is generated by the element
$\ug:=(t_2-t_3\1)\cg_2''$.
\endroster
It is immediate that
\[
\Z[G_{23}]\cg_1\oplus\Z[G_{13}]\cg_2\oplus\Z[G_{12}]\cg_3
 =\joinrel=\Bm_4/\Bm_3,
\label{epi.Bm4}
\]
see~\eqref{rel.bb}; the other relations do not affect this quotient.
Furthermore,
as obviously $\cpp_m(t_3)\cg_1'=\cpp_m(t_3)\cg_2'=\cpp_m(t_1)\cg_3'=0$,
we have an epimorphism
\[
(\Z[G_{01}]/\cpp_m)\cg_1'\oplus(\Z[G_{02}]/\cpp_m)\cg_2'\oplus
 (\Z[G_{03}]/\cpp_m)\cg_3'\longonto\Bm_3/\Bm_2.
\label{epi.Bm3}
\]
In $\Bm_2$, we have a relation
\[*
\cg_1''=\cg_2''+\cg_3'';
\]
it is the linear combination
$(t_1-1)\eqref{rel.a1}-(t_2-1)\eqref{rel.a2}-(t_3-1)\eqref{rel.a3}$.
Multiplying this by $(t_2-t_3\1)$ and using~\eqref{rel.bb}, we have
\[*
\ug:=(t_2-t_3\1)\cg_2''=-(t_2-t_3\1)\cg_3''.
\]
Hence, using~\eqref{rel.bb} again, we obtain epimorphisms
\begin{gather}
(\Z[G_{3}]/\cpp_m)\cg_2''\oplus(\Z[G_{2}]/\cpp_m)\cg_3''\longonto\Bm_2/\Bm_1,
 \label{epi.Bm2}\\
(\Z[G_{1}]/\cpp_m)\ug\longonto\Bm_1\quad\text{(for $m$ odd)}.
 \label{epi.Bm1.odd}
\end{gather}
If $m=2k$ is even, arguing as in~\cite{degt:Fermat} we can
refine~\eqref{epi.Bm1.odd} to
\[
(\Z[G_{1}]/\cpp_k(t^2))\ug\longonto\Bm_1\quad\text{(for $m=2k$ even)},
 \label{epi.Bm1.even}
\]
where $t:=t_0=t_1=t_2\1=t_3\1$. Indeed,
since $t_2\ug=t_3\ug=t_1\1\ug$,
by induction for $r\in\Z$ we have
\[*
t_2^r\cg_2''=t_1^r\cg_2''+t_2^{1-r}\cpp_r(t_2^2)u,\quad
t_3^r\cg_3''=t_1^r\cg_3''-t_3^{1-r}\cpp_r(t_3^2)u.
\]
Summing up and using the fact that $\cpp_m(t_1)\cg_2''=\cpp_m(t_2)\cg_2''=0$
and the identity
\[*
t^{m-2}\sum_{r=0}^{m-1}t^{1-r}\cpp_r(t^2)=
 t\cpp_{k-1}(t^2)\cpp_m(t)+\cpp_k(t^2),\quad m=2k
\]
(which is easily established by multiplying both sides by $t^2-1$),
we immediately conclude that $\cpp_k(t_2^2)\ug=0$.

Since $\Ga=m\in\N_+$, we have isomorphisms
$G_{ij}\cong G_i\cong\Z/m$ and, hence, all rings
$\Z[G_*]/\cpp_m$ in~\eqref{epi.Bm2} and~\eqref{epi.Bm1.odd}
are free abelian groups of rank $m-1$.
If $m=2k$ is even, the ring $\Z[G_{1}]/\cpp_k(t^2)$ in~\eqref{epi.Bm1.even}
is a free abelian group of rank $m-2$.
Thus,
summing up, we have
$\ell(\Bm\of\Ga)\le9m-6-\Gd\of\Ga$.
On the other hand, due to~\eqref{eq.rank} and~\cite{Shioda:Fermat},
$\rank\Bm\of\Ga=9m-6-\Gd\of\Ga$.
Hence,
$\Tors\Bm\of\Ga=0$.
\qed

\corollary[of the proof]\label{cor.simple.relation}
The $\GL\ofm$-module~$\Bm\ofm$ can be defined by relations
\eqref{rel.bb} and
$\cg_1''=\cg_2''+\cg_3''$, where $\cg_i''$ are the elements introduced
in \autoref{proof.Fermat}.
Furthermore, all epimorphisms
in~\eqref{epi.Bm4}--\eqref{epi.Bm1.even} are isomorphisms.
\done
\endcorollary

\remark\label{rem.simple.relation}
\autoref{cor.simple.relation} does not extend to other finite
quotients, \CF. \autoref{s.cyclic}.
\endremark

\subsection{Proof of \autoref{th.main}}\label{proof.main}
In view of~\eqref{eq.rank}, the rank $\rank\KF\Ga$ can be computed as
$\dim_\C(\Bm\of\Ga\otimes\C)+3$.
The group algebra $\C[\GG/m\GG]$ is semisimple, and we have
\[*
\Bm\ofm\otimes\C=
\Bm_1\otimes\C\oplus(\Bm_2/\Bm_1)\otimes\C
 \oplus(\Bm_3/\Bm_2)\otimes\C\oplus(\Bm_4/\Bm_3)\otimes\C,
\]
see \autoref{proof.Fermat}.
The rank formula in the theorem
is obtained by tensoring this expression by $\C[G]$ and
using
isomorphisms~\eqref{epi.Bm4}--\eqref{epi.Bm1.even}.

Let $(i,j,k)$ be a permutation of $(1,2,3)$, and introduce the following
parameters, measuring the `inhomogeneity' of $\Ker\Ga$:
\roster*
\item
$m_i$ is the order of the image~$\Ga(t_i)$ in~$G$;
\item
$n_i$ is the order of
the image of $t_i$ (or $t_0$) in $G/\Ga(t_0t_i)=G/\Ga(t_jt_k)$;
\item
$n_{jk}$ is the order of
the image of $t_j$ (or $t_k$) in $G/\Ga(t_0t_i)=G/\Ga(t_jt_k)$;
\item
$\bar n_i:=n_i/\ls|G_{jk}|=n_{jk}/\ls|G_{0i}|$;
\item
$p_i:=\gcd(n_i,n_{jk})$ and $\bar p_i:=p_i/\ls|G_l|$, $i=2,3$, $i+l=5$;
\item
$\bar q:=\gcd(p_2,p_3)/\ls|G_1|$.
\endroster
It is not difficult to see that
all $\bar n_i$, $\bar p_i$, and $\bar q$ are integers.
If
$\Gd\of\Ga=1$,
introduce also
\roster*
\item
$\bar s:=s/\ls|G_1|$, where
$s:=\gcd(s_2,s_3)$ and $s_i:=\lcm(n_i,m_i)$, $i=2,3$.
\endroster
Note that $\bar s$ is an integer and $\bar q\divides|\bar s$.
If $\Gd\of\Ga=0$, we merely let $\bar s:=1$.

\lemma\label{lem.torsion}
There is a filtration
$0=T_0\subset T_1\subset T_2\subset T_3:=\Tors\Bm\of\Ga$
such that the quotient
groups $T_i/T_{i-1}$, $i=1,2,3$, are subquotients of
\[*
(\Z/\bar q)\oplus(\Z/\bar s),\quad
(\Z/\bar p_2)\oplus(\Z/\bar p_3),\quad
(\Z/\bar n_1)\oplus(\Z/\bar n_2)\oplus(\Z/\bar n_3),
\]
respectively.
In particular, $\ell(\Tors\Bm\of\Ga)\le6+\Gd\of\Ga$.
\endlemma

\proof
Over $\GL\ofm$, the tensor product does not need to be exact, but we still
have an epimorphism $\Bm\ofm\otimes_{\GL\ofm}\GL\of\Ga\onto\Bm\of\Ga$, which
induces an epimorphism of the torsion groups (as the ranks of the two
modules, regarded as abelian groups, are equal).
Using the same filtration as in \autoref{proof.Fermat}, we obtain
epimorphisms \eqref{epi.Bm4}--\eqref{epi.Bm1.even},
which
also induce epimorphisms of the torsion subgroups.
Then, define the member $T_i\subset\Tors\Bm\of\Ga$ of the filtration
as the image of $\Bm_i$, $i=0,1,2,3$.

The group rings $\Z[G_*]$ in~\eqref{epi.Bm4} are torsion free; hence, indeed,
$T_3=\Tors\Bm\of\Ga$.
Let $(i,j,k)$ be a permutation of $(1,2,3)$.
In~\eqref{epi.Bm3}, each generator~$\cg_i'$ is
annihilated by $\cpp_{n_{jk}}(t_j)$,
and we can refine the corresponding summand to
$(\Z[G_{0i}]/\cpp_{n_{jk}})\cg_i'$.
Let $r_i:=\ls|G_{0i}|$ be the order of the cyclic group~$G_{0i}$.
Then $\cpp_{n_{jk}}=\bar n_i\cpp_{r_i}$ in $\Z[G_{0i}]$, and
$\Z[G_{0i}]/\cpp_{r_i}$ is a free abelian group of rank $r_i-1$.
Hence, $\Tors(\Z[G_{0i}]/\cpp_{n_{jk}})\cg_i'$
is a cyclic group $\Z/\bar n_i$;
more precisely,
\[*
\ord\bigl((t_k^{r_i}-1)\cg_i\bigr)\
 \text{in $\Bm\of\Ga/\Bm_2$ divides $\bar n_i$, where $r_i:=\ls|G_{0i}|$}.
\]
Tensoring this element by~$\C$, one can see that it does have finite
order in $\Bm\of\Ga/\Bm_2$ but, in general, not in $\Bm\of\Ga$.

A similar argument applies to~\eqref{epi.Bm2} and~\eqref{epi.Bm1.odd}.
In~\eqref{epi.Bm2}, the summand generated by~$\cg_i''$ is refined to
$(\Z[G_l]/\cpp_{p_i})\cg_i''$, $l:=5-i$,
the torsion of which is $\Z/\bar p_i$:
\[*
\ord\bigl((t_i-1)(t_l^{r_l}-1)\cg_i\bigr)\
 \text{in $\Bm\of\Ga/\Bm_1$ divides $\bar p_i$, $i=2,3$, where $r_l:=\ls|G_l|$}.
\]
In~\eqref{epi.Bm1.odd}, the module refines to
$(\Z[G_1]/\cpp_q)\ug$, and we have
\[*
\ord\bigl((t_2-t_1)(t_2^{r}-1)(t_3-1)\cg_2\bigr)\
 \text{in $\Bm\of\Ga$ divides $\bar q$, where $r:=\ls|G_1|$}.
\]

If
$\Gd\of\Ga=1$ (equivalently, if both $m=2k$ and $\ls|G_1|=2l$ are even),
we use~\eqref{epi.Bm1.even} instead
of~\eqref{epi.Bm1.odd}. In addition to $\cpp_q(t)\ug=0$, we also have
$\cpp_{s/2}(t^2)\ug=0$, \CF. the end of \autoref{proof.Fermat}.
Since $\cpp_{s/2}(t^2)=\bar s\cpp_l(t^2)$ and
$\cpp_q(t)=\bar q(t+1)\cpp_l(t^2)$ in $\Z[G_1]$,
we obtain an extra torsion term:
\[*
\ord\bigl((t_2-t_1)\cpp_r(-t_2)(t_3-1)\cg_2\bigr)\
 \text{in $\Bm\of\Ga$ divides $\bar s$, where $r:=\ls|G_1|$}.
\]
Comparing the ranks, we conclude that the elements indicated above exhaust
all torsion that may be present in $\Bm\of\Ga$.
\endproof

\remark\label{rem.torsion}
Note that \autoref{lem.torsion} is merely an estimate on the size of
$\TF\Ga$. In particular, its conclusion depends on the order of the
indices, and one may get a better estimate by permuting the indices
$(0,1,2,3)$, \CF. \autoref{rem.non-split} and \autoref{s.cyclic}.
\endremark

Denote by $\Gs\:\tFermat\ofm\to\Fermat\ofm$ the Fermat surface~$\Fermat\ofm$
blown up so
that the projection $\pi\:\tFermat\ofm\to\Fermat\of\Ga$ is regular,
and let $\tdv\ofm:=\Gs^*\dv\ofm$.

\lemma\label{lem.pi}
The maps
\[*
\NS(\Fermat\of\Ga)\overset{\smash{\pi^*}}\longto
 \NS(\tFermat\ofm)\overset{\smash{\pi_*}}\longto\NS(\Fermat\of\Ga)
\]
respect the subgroups
$\SF\Ga\subset\NS(\Fermat\of\Ga)$ and
$\SS{\tdv\ofm}\subset\NS(\tFermat\ofm)$.
The composition
$\pi_*\circ\pi^*\:\NS(\Fermat\of\Ga)\to\NS(\Fermat\of\Ga)$
is the multiplication by $d:=m^3\!/\ls|G|$.
\endlemma

\proof
The first statement is immediate from the definition of the divisors
involved: set-theoretically, one has $\dv\of\Ga=\pi(\tdv\ofm)$ and
$\tdv\ofm=\pi\1(\dv\of\Ga)$. The second statement is well known: since
$\pi$ is a generically finite-to-one map of degree~$d$, the assertion is
geometrically obvious for the class of an irreducible curve
$C\subset\Fermat\of\Ga$ not contained in the ramification locus; then, it
remains to observe that $\NS(\Fermat\of\Ga)$ is generated by such classes
(\eg, very ample divisors).
\endproof

By \autoref{lem.pi}, we have induced maps
\[*
\NS(\Fermat\of\Ga)/\SF\Ga\overset{\smash{\pi^*}}\longto
 \NS(\tFermat\ofm)/\SS{\tdv\ofm}\overset{\smash{\pi_*}}\longto
 \NS(\Fermat\of\Ga)/\SF\Ga
\]
whose composition $\pi_*\circ\pi^*$
is the multiplication by~$d$. Since the group in the middle
is torsion free, see \autoref{th.free}\iref{free.Fermat}
and~\eqref{eq.sigma},
the group $\TF\Ga\subset\Ker\pi^*$ is annihilated by~$d$. Together with the
estimate on $\ell(\TF\Ga)$ given by \autoref{lem.torsion}, this completes the
proof of \autoref{th.main}.
\qed

\subsection{Proof of \autoref{th.free}\iref{free.diagonal}}\label{proof.diagonal}
The statement follows from \autoref{lem.torsion}, as one obviously has
$\bar n_i=\bar p_i=\bar q=1$, $i=1,2,3$, and
\[*
\bar s=\lcm_{1\le i<j\le3}(\gcd(m_i,m_j))/\gcd(m_1,m_2,m_3).
\]
In fact, using \autoref{cor.simple.relation}, one can easily show that
$\Tors(\Bm\ofm\otimes_\GL\GL\of\Ga)=\Z/\bar s$.
Furthermore,
numeric examples suggest that
$\Bm\ofm\otimes_\GL\GL\of\Ga=\Bm\of\Ga$
in the diagonal case, see \autoref{s.diagonal}.
However, I do not know a proof of the latter statement.
\qed

\subsection{Proof of \autoref{ad.lattice}}\label{proof.lattice}
It suffices to show that $\SF\Ga$
is a nondegenerate lattice; then,
in addition to the obvious inclusion
$\KF\Ga\subset\ker H_2(\dv\of\Ga)$,
we would also have the converse statement
$\KF\Ga\supset\ker H_2(\dv\of\Ga)$.

Consider the subspace $\SF\Ga\otimes\Q\subset H_2(\Fermat\of\Ga;\C)$ and
recall that each multi-eigenspace (isotypical component of the $G$-action) in
$H_2(\Fermat\of\Ga;\C)$ has dimension at most~$1$,
see~\cite{Shioda:Fermat,Shioda:Delsarte}. Hence, $\SF\Ga\otimes\Q$ is a
direct sum of whole eigenspaces,
which
are obviously nondegenerate and orthogonal.
\qed

\subsection{Proof of \autoref{cor.generated}}\label{proof.generated}
According to~\cite{Shioda:Fermat}, for any integer $m\in\N_+$
prime to~$6$, one has $\NS(\Fermat\ofm)\otimes\Q=\SF{m}\otimes\Q$. Then, by
\autoref{lem.pi}, a similar identity
$\NS(\Fermat\of\Ga)\otimes\Q=\SF\Ga\otimes\Q$
holds for any finite quotient $\Ga\:\GG\onto G$ with $\ls|G|$ prime to~$6$.
It remains to observe that, for each surface $\Fermat\of\Ga$ as in the
statement,
\roster*
\item
$\pi_1(\Fermat\of\Ga)=0$, see \autoref{th.pi1}; hence,
$\Pic\Fermat\of\Ga=\NS(\Fermat\of\Ga)$, and
\item
$\TF\Ga=0$, see \autoref{th.free}.
\endroster
(If $\Ga$ is diagonal, \CF. \autoref{th.free}\iref{free.diagonal}, the
assumption that $\ls|G|$ is prime to~$6$ implies also that $\Gd\of\Ga=0$.)
The last statement follows from \autoref{ad.lattice}.
\qed

\section{Cyclic Delsarte surfaces}\label{S.cyclic}

Throughout this section, we fix an epimorphism $\Ga\:\GG\onto G$ and assume
that $G$ is a finite cyclic group,
$\ls|G|=m$.

\subsection{The setup}\label{s.setup}
Fix a generator~$t$ of~$G$ and let $\Ga(t_i)=t^{m_i}$, $i=0,1,2,3$.
Strictly speaking, $m_0,m_1,m_2,m_3$ are elements of $\Z/m$, but
it is more convenient to regard
them as nonnegative integers.
Then $m_0+m_1+m_2+m_3=0\bmod m$ and
\[
\gcd(m,m_1,m_2,m_3)=1.
\label{eq.gcd}
\]
For $i\ne j$, let $m_{ij}:=\gcd(m,m_i+m_j)$.
We have $m_{ij}=m_{kl}$ whenever $(i,j,k,l)$ is
a permutation of $(0,1,2,3)$, \ie, there are
three essentially distinct parameters $m_{ij}$.

It is easy to see that $\Gd\of\Ga=1$ if and only if
$m=0\bmod2$ and $m_1m_2m_3=1\bmod2$.
In view of~\eqref{eq.gcd},
\[
\gcd(m_{12},m_{13},m_{23})=2^{\Gda}.
\label{eq.gcd.ij}
\]
The following statement is an immediate consequence of~\eqref{eq.gcd}
and~\eqref{eq.gcd.ij}.

\lemma\label{lem.exclusive}
For a divisor $d\divides|m$, $d>2$, the following two conditions
\roster
\item\label{cyclic.1}
$d\divides|m_i$ and $d\divides|m_j$ for some $0\le i<j\le3$, or
\item\label{cyclic.2}
$d\divides|m_{ij}$ and $d\divides|m_{ik}$ for some permutation $(i,j,k)$ of
$(1,2,3)$,
\endroster
are mutually exclusive.
Furthermore, $d$ may satisfy
either~\iref{cyclic.1} for at most one pair $i<j$ or~\iref{cyclic.2}
for at most one value of $i\in\{1,2,3\}$.
\done
\endlemma

\subsection{Proof of \autoref{th.pi1} for cyclic Delsarte surfaces}\label{proof.pi1.cyclic}
Due to the general expression for $\pi_1(\Fermat\of\Ga)$ given by
\autoref{th.pi1},
it suffices to show that, in the ring $\Z/m$, each
solution to the equation $r_1m_1+r_2m_2+r_3m_3=0$ decomposes into a sum of
solutions with at least one unknown $r_i=0$.
Since $\Z/m=\bigoplus_q\Z/q$, the summation running over all maximal prime
powers $q\divides|m$, we can assume that $m$ itself is a prime power. Then,
due to~\eqref{eq.gcd}, at least one coefficient~$m_i$ is prime to~$m$.
If, for example,
$\gcd(m,m_1)=1$, \ie, $m_1$ is invertible in~$\Z/m$, we obtain an
equivalent equation $r_1=-r_2n_2-r_3n_3$, where $n_i:=m_im_1\1$, $i=1,2$, for
which the decomposition statement is obvious.
\qed

\subsection{Invariant factors}\label{s.inv}
In the rest of this section,
we prove \autoref{th.free}\iref{free.cyclic}
by analyzing the structure of the module
$\Am\of\Ga$ (see \autoref{rem.cyclic} for an explanation).
Introduce the notation
\[*
\cp:=t^m-1,\quad
\cp_i:=t^{m_i}-1,\quad
\cp_{ij}:=t^{m_{ij}}-1,\quad
i,j=0,1,2,3,\quad i\ne j.
\]
Recall that, for $p,q\in\Z$, one has $\gcd(t^p-1,t^q-1)=t^{\gcd(p,q)}-1$.
Hence, the polynomials introduced
are subject to the following divisibility relations:
\[
\gathered
\cp_{ij}\divides|\cp\quad\text{for all $i\ne j$ (by the definition of $m_{ij}$)},\\
\gcd(\cp,\cp_1,\cp_2,\cp_3)=\co:=t-1\quad\text{(see~\eqref{eq.gcd})},\\
\gcd(\Gs_i,\Gs_j,\Gs_{ik})=\gcd(\Gs_i,\Gs_{ij},\Gs_{ik})=\co\quad
 \text{for $\{i,j,k\}=\{1,2,3\}$},\\
\gcd(\cp_{12},\cp_{13},\cp_{23})=\cd\co,\quad\cd:=(t+1)^{\Gda}\quad
 \text{(see~\eqref{eq.gcd.ij})}.
\endgathered
\label{eq.divisibility}
\]
(The third relation follows from the similar relations for the
exponents~$m_*$, which, in turn, are consequences of~\eqref{eq.gcd}.)
These relations hold in the following \emph{ideal sense}:
the ideal generated in
$\tGL:=\Z[t^{\pm1}]$ by the polynomials in
$\gcd(\ldots)$
the left hand side of a relation
equals the
ideal generated by the polynomial in the right hand side.
In particular, they hold over~$\Z$ as well as over any field.

We regard $\Am\of\Ga$ as an $\tGL$-module. It is generated by
$\ag_1$, $\ag_2$, $\ag_3$, $\cg_1$, $\cg_2$, $\cg_3$, and the defining
relations are \eqref{rel.a1}--\eqref{rel.a3}
with $t_i=t^{m_i}$, $i=1,2,3$, and
\[*
\cp\ag_1=\cp\ag_2=\cp\ag_3=\cp_{23}\cg_1=\cp_{13}\cg_2=\cp_{12}\cg_3=0.
\]
(The first three relations make $\Am\of\Ga$ a $\Z[G]$-module, and the last
three are~\eqref{rel.bb}
combined with
$\cp\cg_i=0$, $i=1,2,3$.)
The relations in~$\Am\of\Ga$ are represented by the matrix
\[\label{eq.M}
M:=\left[\begin{matrix}
0&\cp_3&-\cp_2&\cp_3&0&0\\
\cp_3&0&-\cp_1&0&\cp_3&0\\
-\cp_2&\cp_1&0&0&0&\cp_1\\
\cp&0&0&0&0&0\\
0&\cp&0&0&0&0\\
0&0&\cp&0&0&0\\
0&0&0&\cp_{23}&0&0\\
0&0&0&0&\cp_{13}&0\\
0&0&0&0&0&\cp_{12}
\end{matrix}\right].
\]

Given a field~$\kk$,
the reduction $\Am\of\Ga\otimes\kk$
is a finitely generated module over the principal ideal domain
$\kk\tGL:=\Z[t^{\pm1}]\otimes\kk=\kk[t^{\pm1}]$;
hence, it decomposes into direct sum of cyclic modules,
\[*
\Am\of\Ga\otimes\kk\cong
\kk\tGL/\inf_1\oplus
\ldots
\oplus\kk\tGL/\inf_6,
\]
where $\inf_1,\ldots,\inf_6$ are the \emph{invariant factors}
of~$M\otimes\kk$, \ie,
the diagonal elements of the Smith normal form of the matrix.
Recall that
$\inf_1\divides|\inf_2\divides|\inf_3\divides|\inf_4\divides|\inf_5\divides|\inf_6$
are elements of~$\kk\tGL$ that can be found as
$\inf_r=(\gcd S_r)/(\gcd S_{r-1})$, $r=1,\ldots,6$, where $S_r$ is the set of
all $(r\times r)$-minors of~$M\otimes\kk$.


All nontrivial minors of~$M$ are products of polynomials of the form
$(t^s-1)$. Computing all $(r\times r)$-minors, $r=1,\ldots,6$, we obtain six
lengthy sequences~$S_r$. Since we are interested in the greatest common divisors
only, we use~\eqref{eq.divisibility}
(in the ideal sense as explained above)
and simplify these sequences as described below.

Whenever a sequence~$S$ contains a subsequence of the form
\roster*
\item
$\Gb\cp$, $\Gb\cp_1$, $\Gb\cp_2$, $\Gb\cp_3$, or
\item
$\Gb\cp_i$, $\Gb\cp_j$, $\Gb\cp_{ik}$ for some $\{i,j,k\}=\{1,2,3\}$, or
\item
$\Gb\cp_i$, $\Gb\cp_{ij}$, $\Gb\cp_{ik}$ for some $\{i,j,k\}=\{1,2,3\}$,
\endroster
where $\Gb$ is a common factor, one can append to~$S$ the product
$\Gb\co$.
After all such additions have been made, one can shorten~$S$ by removing all
nontrivial
multiples of any element $\Gb\in S$.
We repeat
these two steps
until $S$ stabilizes, and then apply a similar
procedure, replacing each subsequence
$\Gb\cp_{12}$, $\Gb\cp_{13}$, $\Gb\cp_{23}$ with
the product $\Gb\cd\co$.
Denoting by~$S_r'$ the result of the simplcication, we have
\[
\gathered
S_1'=\{\co\},\quad
S_2'=\{\co^2\},\quad
S_3'=\{\co^3\},\quad
S_4'=\{\cd\co^4\},\\
S_5'=\{\cp\cd\co^4,\cp_{12}\cp_{13}\cp_{23}\co^2,
 \cp_2\cp_3\cp_{12}\cp_{23}\co,\cp_1\cp_3\cp_{12}\cp_{13}\co,
 \cp_1\cp_3\cp_{13}\cp_{23}\co\}.
\endgathered
\label{eq.minors}
\]
Another observation is the fact that $S_6$ is a subset of
$\{\cp\Gb\,|\,\Gb\in S_5\}$; hence, one has
$\cp(\gcd S_5)\divides|\gcd S_6$. On the other hand, $\Am\of\Ga$ is a
$\Z[G]$-module and all its invariant factors are divisors of~$\cp$.
Taking into account~\eqref{eq.minors}, we easily obtain all invariant
factors (in any characteristic) except~$\inf_5$:
\[
\inf_1=\inf_2=\inf_3=\co,\quad
\inf_4=\cd\co,\quad
\inf_6=\cp.
\label{eq.factors}
\]

\subsection{The factor~$\inf\sb5$: the case \pdfstr{k=Q}{$\kk=\Q$}}\label{s.p=0}
Let $\bcp_*:=\cp_*/\co$ and cancel the common factor~$\co^5$,
converting~$S_5'$ to the union
\[*
S_5'':=\{\bcp\cd\}\cup\{\bcp_{12}\bcp_{13}\bcp_{23},
 \bcp_2\bcp_3\bcp_{12}\bcp_{23},\bcp_1\bcp_3\bcp_{12}\bcp_{13},
 \bcp_1\bcp_3\bcp_{13}\bcp_{23}\}.
\]
Over $\Q$,
the irreducible factors of~$\cp$
are distinct cyclotomic polynomials $\cf_d$, $d\divides|m$, and a
factor $\cf_d$,
$d>2$, may appear in $\gcd S_5''$
at most once.
Since $\cf_d\divides|\bcp_{12}\bcp_{13}\bcp_{23}$, one has
$d\divides|m_{ij}$ for some $1\le i<j\le3$.
It remains to consider the three possibilities case-by-case and analyze the
remaining three elements of~$S_5''$. Using the relations between~$m_*$
(mainly, the fact that $\gcd(m_i,m_{ij})=\gcd(m,m_i,m_j)$), we arrive at the
following restrictions to~$d$:
\roster
\item\label{p=0.1}
$d\divides|m_i$ and $d\divides|m_j$ for some $1\le i<j\le3$, or
\item\label{p=0.2}
$d\divides|m_{ij}$ and $d\divides|m_{ik}$ for some permutation $(i,j,k)$ of
$(1,2,3)$, or
\item\label{p=0.3}
$d\divides|m_i$ and $d\divides|m_{jk}$ for some permutation $(i,j,k)$ of
$(1,2,3)$.
\endroster
The substitution $m_{jk}\mapsto m_{i0}=-m_{jk}\bmod m$ converts~\iref{p=0.3}
to~\iref{p=0.1} with $(i,j)=(i,0)$.
Hence, $\gcd S_5'=\inf_5\cd\co^4$ with
\[\label{eq.f5}
\inf_5=\prod\cf_d(t),
\]
where
the product runs over all divisors $d\divides|m$ satisfying
conditions~\iref{cyclic.1} or~\iref{cyclic.2} in \autoref{lem.exclusive}.
(In the special case $d=2$ and $\Gda=1$, the greatest common divisor
contains two copies of $(t+1)$; one of them is~$\cd$, and the other is
counted in the product~\eqref{eq.f5} for~$\inf_5$.
An extra factor $\co=\cf_1(t)$ is also counted in the product.)

\remark\label{rem.cyclic.rank}
According to~\eqref{eq.factors} and~\eqref{eq.f5},
$\rank\Am\of\Ga=m+4+\Gda+\sum_d\Euler(d)$,
where $\Euler(d)=\deg\cf_d$ is Euler's totient
function and
the summation runs over all divisors $d\divides|m$
satisfying conditions~\iref{cyclic.1} or \iref{cyclic.2} in
\autoref{lem.exclusive}.
Since $n=\sum_{d\divides|n}\Euler(d)$ for $n\in\N_+$,
this expression translates to
$\rank\Am\of\Ga=m-4-\Gda+\sum_{i<j}d_{ij}+\sum_{i}d_i$
(using \autoref{lem.exclusive} again), where
\roster*
\item
$d_{ij}:=\gcd(m,m_i,m_j)=\ls|G_{ij}|$ for $0\le i< j\le3$, and
\item
$d_i:=\gcd(m_{ij},m_{ik})=\ls|G_i|$ for $i=1,2,3$ and $\{i,j,k\}=\{1,2,3\}$.
\endroster
This agrees with~\eqref{eq.rank} and \autoref{th.main}.
\endremark

\subsection{The factor~$\inf\sb5$: the case \pdfstr{k=F\_5}{$\kk=\FF_p$}}\label{s.p>0}
Fix a prime $p>0$ and compute~$\inf_5$ over~$\FF_p$.
This time, the cyclotomic polynomials~$\cf_d$ may be reducible.
However, for any pair $n,d\in\N_+$ with $\gcd(d,p)=1$ one still has
$\cf_d\divides|(t^n-1)$ if $d\divides|n$ and $\gcd(\cf_d,t^n-1)=1$ otherwise.
Thus, if $p$ is prime to~$m$ (and hence $\cp_m$ is square free),
the computation runs exactly as in \autoref{s.p=0}
and we arrive at~\eqref{eq.f5}.

In general,
let $m_*=m_*'q_*$, where $q_*$ is a power of~$p$ and $m'_*$ is prime to~$p$.
Then,
\smash{$\cp_*=(\cp'_*)^{q_*}$},
where \smash{$\cp'_*:=t^{m'_*}-1$} is square free.
To
reduce the number of cases and
simplify the argument, note that the isomorphism class of the module
$\Am\of\Ga\otimes\FF_p$ and, hence, its invariant factors depend on~$m$
and \emph{unordered} quadruple $(m_0,m_1,m_2,m_3)$ only.
Thus, permuting the indices, we can add to~$S_5'$ all products of the form
$\cp_i\cp_j\cp_{ij}\cp_{ik}\co$, where $(i,j,k)$ runs over all three-element
arrangements of $\{0,1,2,3\}$.
Denote this new set by~$S_5''$.

Let $d'\divides|m'$, $d'>2$. Arguing as in \autoref{s.p=0}, we conclude that
$\cf_{d'}$ divides $\gcd S_5'$ if and only if
\roster
\item\label{p>0.1}
$d'\divides|m'_i$ and $d'\divides|m'_j$ for some $0\le i<j\le3$, or
\item\label{p>0.2}
$d'\divides|m'_{ij}$ and $d'\divides|m'_{ik}$ for some permutation $(i,j,k)$ of
$(1,2,3)$.
\endroster
As in \autoref{lem.exclusive},
the two conditions are mutually exclusive and $d'$ may satisfy
either~\iref{p>0.1} for exactly one pair $i<j$
or~\iref{p>0.2} for exactly one value of~$i$.

In case~\iref{p>0.1}, assume that $(i,j)=(1,2)$ and $q_1=\min(q_1,q_2)$.
Then $d'$ divides~$m'_1$, $m'_2$, $m'_{12}$ and $m'_{03}$,
and $d'$ does not divide any other of~$m'_k$ or~$m'_{kl}$.
Considering the element $\cp_1\cp_3\cp_{13}\cp_{01}\co\in S''_5$,
we see that the multiplicity of $\cf_{d'}$ in $\gcd S_5'$ is at most (and
hence equal to) $q':=\min(q,q_1)$, \ie, the one given by~\eqref{eq.f5} reduced
modulo~$p$.
Indeed, for $\cf_{d'}$, the product in~\eqref{eq.f5} should be
restricted to the divisors of~$m$ of the form $d=d'p^r$.
By the assumption
$q'=\min(q,q_1,q_2)$,
we have $1\le p^r\le q'$.
Since
\[*
\cf_{d'p^r}=(\cf_{d'})^{p^r-p^{r-1}}\quad\text{for $r\ge1$},
\]
the exponents sum up to~$q'$.

In case~\iref{p>0.2}, assume that $(i,j,k)=(1,2,3)$ and $q_{12}\le q_{13}$.
Then $d'$ divides $m'_{12}$, $m'_{13}$, $m'_{03}$ and $m'_{02}$,
and $d'$ does not divide any other of $m'_l$ or $m'_{ln}$.
Considering the element $\cp_1\cp_2\cp_{12}\cp_{01}\co\in S''_5$, as in the
previous case we conclude that
the multiplicity of $\cf_{d'}$ in $\gcd S_5'$ is at most (and
hence equal to) $q_{12}$, \ie, the one given by~\eqref{eq.f5}.

If $d'=1$, the multiplicity of $\cf_1=\co$ (in addition to the four copies
present in each term automatically) is counted by a similar argument, using
the fact that $d=p$ itself satisfies at most one of the two conditions in
\autoref{lem.exclusive} and with at most one parameter set.
The extra multiplicity is $\min(q,q_i,q_j)$ in case~\iref{cyclic.1} or
$\min(q_{ij},q_{ik})$ in case~\iref{cyclic.2},
\ie, again the one given by~\eqref{eq.f5} (where the product is to be
restricted to the divisors $d\divides|m$ that are powers of~$p$).

As in \autoref{s.p=0}, the case
where
$\Gda=1$ and either $d'=2$ or $p=2$ needs
special attention, taking into account the common divisor~$2$ of
all~$m_{ij}$. We leave details to the reader.

Summarizing, we conclude that, for any prime~$p$, the invariant factor~$\inf_5$
of the $\FF_p\tGL$-module
$\Am\of\Ga\otimes\FF_p$
is merely the $(\bmod\,p)$-reduction of~\eqref{eq.f5}.

\subsection{End of the proof of \autoref{th.free}\iref{free.cyclic}}\label{proof.cyclic.full}
For each field $\kk=\Q$ or~$\FF_p$,
\[*
\dim(\Am\of\Ga\otimes\kk)=\deg f_{1,\kk}+\ldots+\deg f_{6,\kk},
\]
where $f_{r,\kk}\in\kk\tGL$, $r=1,\ldots,6$, are the invariant factors of
$\Am\of\Ga\otimes\kk$.
According to \autoref{s.inv}--\autoref{s.p>0},
each $f_{r,\kk}$ is the
reduction to~$\kk$ of the monic polynomial $f_r\in\Z[t]$
given by~\eqref{eq.factors} or~\eqref{eq.f5}.
Hence, $\dim(\Am\of\Ga\otimes\kk)$ does not depend on~$\kk$.
\qed

\section{Examples}\label{S.examples}

In conclusion, we mention a few numeric examples showing the sharpness of most
estimates stated in \autoref{s.results}.
Most examples result from experiments with random matrices, and it appears
that the presence of a nontrivial torsion in~$\Bm\of\Ga$ is quite common.
The input for the computation is a $(3\times3)$-matrix $\GM$ whose rows are
the coordinates (in the basis $t_1,t_2,t_3\in\GG$) of three vectors
generating $\Ker\Ga$.
Usually, this matrix is in the form $\diag(m_1,m_2,m_3)\GM'$, where $\diag$
is a diagonal matrix and $\GM'$ is unimodular:
in the experiments, the diagonal part was fixed while $M'$ was chosen
randomly.

To shorten the display, we represent the isomorphism class
of the finite group $\TF\Ga$ by the vector~$\sm=[a_i]$ of its invariant
factors, so that $\TF\Ga=\bigoplus_i\Z/a_i$.

\subsection{Torsion groups of maximal length}\label{s.max.torsion}
For the finite quotients $\Ga_i$ defined by the matrices $\GM_i:=D\GM'_i$,
where $D:=\diag(1,8,8)$, one has:
\begin{alignat*}3
M'_1&=\mat[ [ 4, 7, 1 ], [ 1, 0, 0 ], [ 0, 1, 0 ] ]:&\quad
 \pi_1(\Fermat\of{\Ga_1})&=\Z/2,&\quad
 \sm&=[ 2, 2, 2, 2, 2, 2, 4 ],\\
M'_2&=\mat[ [ 0, 3, 1 ], [ 1, 0, 0 ], [ 0, 1, 0 ] ]:&\quad
 \pi_1(\Fermat\of{\Ga_2})&=0,&\quad
 \sm&=[ 2, 2, 2, 4 ].
\end{alignat*}
If $D=\diag(1,8,16)$, then
\begin{alignat*}3
M'_3&=\mat[ [ 4, 1, -1 ], [ 1, 1, 0 ], [ 1, 0, 0 ] ]:&\quad
 \pi_1(\Fermat\of{\Ga_3})&=\Z/2,&\quad
 \sm&=[ 2, 2, 2, 4, 4, 4, 4 ],\\
M'_4&=\mat[ [ 6, 1, 2 ], [ 1, 0, 1 ], [ 0, 0, 1 ] ]:&\quad
 \pi_1(\Fermat\of{\Ga_4})&=\Z/4,&\quad
 \sm&=[ 2, 4, 4, 4, 4, 8 ],\\
M'_5&=\mat[ [ 1, 0, 3 ], [ 0, 1, 1 ], [ 0, 0, 1 ] ]:&\quad
 \pi_1(\Fermat\of{\Ga_5})&=0,&\quad
 \sm&=[ 4, 4, 4, 4 ].
\end{alignat*}
If $D=\diag(1,9,9)$ (and hence $\Gd\of\Ga=0$), then
\begin{alignat*}3
M'_6&=\mat[ [ -3, 1, 2 ], [ 1, 0, 0 ], [ 0, 0, 1 ] ]:&\quad
 \pi_1(\Fermat\of{\Ga_6})&=\Z/3,&\quad
 \sm&=[ 3, 3, 3, 3, 3, 9 ],\\
M'_7&=\mat[ [ -1, 1, 1 ], [ 0, 1, 1 ], [ 0, 0, 1 ] ]:&\quad
 \pi_1(\Fermat\of{\Ga_7})&=0,&\quad
 \sm&=[ 3, 3, 9 ].
\end{alignat*}
Finally, for $D=\diag(2,9,9)$ one has
\begin{alignat*}3
M'_8&=\mat[ [ -4, 2, 1 ], [ -3, 1, 0 ], [ 1, 0, 1 ] ]:&\quad
 \pi_1(\Fermat\of{\Ga_8})&=\Z/3,&\quad
 \sm&=[ 3, 3, 3, 3, 3, 3, 9 ],\\
M'_9&=\mat[ [ 3, 2, 0 ], [ 1, 1, 0 ], [ 3, 0, -1 ] ]:&\quad
 \pi_1(\Fermat\of{\Ga_9})&=0,&\quad
 \sm&=[ 3, 3, 3, 9 ].
\end{alignat*}

\remark\label{rem.non-split}
In most examples considered in this section, the estimate given by
\autoref{lem.torsion} does depend on the order of the indices, \CF.
\autoref{rem.torsion}; often, even the best bound is larger than the actual
size $\ls|\TF\Ga|$.
In many cases, the epimorphism
$\Bm\ofm\otimes_\GL\GL\of\Ga\onto\Bm\of\Ga$ is not an isomorphism,
\CF. \autoref{rem.simple.relation}.
Note also that, for the finite quotient~$\Ga_4$, one has
\[*
T_2/T_0\cong T_3/T_2\cong\Z/4\oplus\Z/4\oplus\Z/4,
\]
\CF. \autoref{lem.torsion}, whereas $\exp\TF{\Ga_4}=8$.
\endremark

\subsection{The case of~$\ls|G|$ prime to~$6$}\label{s.prime.6}
In this case, one always has $\Gd\of\Ga=0$.
Let $\Ga_i$ be defined by a matrix $\GM_i:=D\GM'_i$.
If $D=\diag(1,5,25)$, then one has:
\begin{alignat*}3
M'_1&=\mat[ [ 2, -1, 6 ], [ 1, 0, 1 ], [ 0, 0, 1 ] ]:&\quad
 \pi_1(\Fermat\of{\Ga_1})&=\Z/5,&\quad
 \sm&=[ 5, 5, 5, 5, 5, 5 ],\\
M'_2&=\mat[ [ 2, 0, -1 ], [ 4, 1, -1 ], [ 1, 0, 0 ] ]:&\quad
 \pi_1(\Fermat\of{\Ga_2})&=0,&\quad
 \sm&=[ 5, 5, 5 ].
\end{alignat*}
If $D=\diag(1,7,7)$, then
\begin{alignat*}3
M'_3&=\mat[ [ 1, 2, 5 ], [ 0, 0, 1 ], [ 1, 1, 0 ] ]:&\quad
 \pi_1(\Fermat\of{\Ga_3})&=\Z/7,&\quad
 \sm&=[ 7, 7, 7, 7, 7, 7 ],\\
M'_4&=\mat[ [ 1, 0, 2 ], [ 1, 0, 1 ], [ 3, 1, 0 ] ]:&\quad
 \pi_1(\Fermat\of{\Ga_4})&=0,&\quad
 \sm&=[ 7, 7, 7 ].
\end{alignat*}

\remark\label{rem.pi1}
Examples found in \autoref{s.max.torsion} and \autoref{s.prime.6} suggest
that, under the additional assumption that $\pi_1(\Fermat\of\Ga)=0$, we have
a better bound $\ell(\TF\Ga)\le3+\Gd\of\Ga$.
It also appears that $\exp\TF\Ga$ divides $\height\Ga$.
I do not know a proof of these facts.
\endremark

\subsection{Diagonal Delsarte surfaces}\label{s.diagonal}
Tested were the diagonal finite quotients
\[*
\Ga=(2,4,4),\ (2,6,6),\ (2,8,8),\ (4,6,12).
\]
In all cases,
the obvious epimorphism
$\Bm\ofm\otimes_\GL\GL\of\Ga\onto\Bm\of\Ga$ is an isomorphism, \ie, the
torsion $\Tors\Bm\of\Ga$ is maximal allowed by
\autoref{th.free}\iref{free.diagonal},
see \autoref{proof.diagonal}.

\subsection{Cyclic Delsarte surfaces}\label{s.cyclic}
The last example illustrates Remarks~\ref{rem.simple.relation}
and~\ref{rem.torsion},
showing that, in general, one may need to deal with the whole module
$\Am\of\Ga$ when computing the torsion.
Let $\Ga\:\GG\onto G$ be the finite quotient defined by the matrix
\[*
M:=\mat[ [  1,  1,  0 ], [  3,  0,  3 ], [  0,  0,  4 ] ].
\]
It is immediate that $m=12$ and $G\cong\Z/m$ is a cyclic group; hence,
$\Tors\Bm\of\Ga=0$, see~\eqref{eq.Tors} and
\autoref{th.free}\iref{free.cyclic}.

Let $\Bm'\of\Ga:=\Bm\ofm\otimes_\GL\GL\of\Ga$; by
\autoref{cor.simple.relation}, this $\GL\of\Ga$-module is defined
by~\eqref{rel.bb} and relation $\cg_1''=\cg_2''+\cg_3''$.
Consider the filtrations $\Bm_i\subset\Bm\of\Ga$ and
$\Bm'_i\subset\Bm'\of\Ga$, $i=0,\ldots,4$,
defined as in \autoref{proof.Fermat}.
Then, a straightforward computation shows that
$\Tors(\Bm_3/\Bm_2)=\Z/4\oplus\Z/2$ whereas
$\Tors(\Bm_3'/\Bm_2')=\Z/4\oplus\Z/4\oplus\Z/2$ (as predicted by
\autoref{lem.torsion}); hence, $\Bm\of\Ga\ne\Bm'\of\Ga$, \CF.
\autoref{rem.simple.relation}.

Furthermore,
$\bar p_2=\bar p_3=2$ and $\bar q=\bar s=1$ and,
in agreement with \autoref{lem.torsion}, we have
$\Tors\Bm_2=\Tors\Bm_2'=\Z/2\oplus\Z/2$.
However, permuting the indices to $(0,2,1,3)$
(\CF. \autoref{rem.torsion}), we obtain a better bound:
this time $\bar p_2=\bar p_3=\bar q=\bar s=1$ and, hence,
$\Tors\Bm_2=\Tors\Bm_2'=0$.


\remark\label{rem.cyclic}
This example explains also why, in the proof of
\autoref{th.free}\iref{free.cyclic} in \autoref{S.cyclic}, we had to consider
the matrix~\eqref{eq.M} with rather long sequences of minors instead of a
much simpler matrix given by \autoref{cor.simple.relation}:
the latter just would not work, as the corresponding module \emph{may} have
torsion.
\endremark

\let\.\DOTaccent
\def\cprime{$'$}
\bibliographystyle{amsplain}
\bibliography{degt}

\end{document}